\documentclass[12pt,reqno]{amsart}
\usepackage{amscd,amssymb,amsmath}
\usepackage[arrow,matrix]{xy}
\usepackage{graphicx}
\usepackage{cite}
\usepackage{ifpdf}
\usepackage{tikz}

\newtheorem{thm}[subsection]{Theorem}

\newtheorem{prop}[subsection]{Proposition}
\newtheorem{cor}[subsection]{Corollary}

\newtheorem{rk}[subsection]{Remark}
\newtheorem{defn}[subsection]{Definition}
\newtheorem{ex}[subsection]{Example}

\numberwithin{equation}{section} 

\def \t {\theta}
\def \w {\omega}
\def \> {\Rightarrow}
\def \0 {\emptyset}

\def \e {\varepsilon}

\def\b {\beta}

\newcommand{\Z}{\mathbb{Z}}
\newcommand{\Q}{\mathbb{Q}}

\newcommand{\N}{\mathbb{N}}
\title [On classification of $p$-adic Leibniz algebras]
{Classification of $p$-adic 6-dimensional filiform Leibniz algebras by solution of $x^q=a$}

\author {M. Ladra \and B.A. Omirov \and U.A. Rozikov}

\thanks{ The first author was supported by Ministerio
de Ciencia e Innovaci\'on (European FEDER support included), grant
MTM2009-14464-C02, and by Xunta de Galicia, grant Incite09 207 215
PR. The second and third authors thank the Department of Algebra, University
of Santiago de Compostela, Spain,  for providing financial support
to their visit to the Department}

\begin{document}

\begin{abstract} In this paper we study the $p$-adic equation $x^q=a$ over the field of $p$-adic numbers.
We construct an algorithm of calculation of criteria of solvability in the case of $q=p^m$ and present
a computer program to compute the criteria for fixed value of $m \leq p-1$. Moreover, using
this solvability criteria for $q=2,3,4,5,6$, we classify $p$-adic 6-dimensional filiform Leibniz algebras.
\end{abstract}


\maketitle

\section{Introduction} \label{sec:intro}

For a given prime $p$, the field $\Q_p$ of $p$-adic numbers is a completion of the rational numbers.
The field $\Q_p$ is also given by a topology derived from a metric, which is itself derived from an alternative
valuation on the rational numbers. This metric space is complete in the sense that every Cauchy sequence converges
to a point in $\Q_p$. This is what allows the development of calculus on $\Q_p$ and it is the interaction
of this analytic and algebraic structure which gives the $p$-adic number systems their power and utility.

 A long time after the discovery of the $p$-adic numbers, they were
mainly considered objects of pure mathematics. However, numerous applications
of these numbers to theoretical physics have been proposed in papers \cite{ADFV,FW,MP} to
quantum mechanics \cite{Kh1}, to $p$-adic - valued physical observables \cite{Kh1} and many
others \cite{Kh2,VVZ}.

In the study of a variety of algebras the classification of algebras in low dimensions plays a crucial role,
since it is helpful to observe some properties of the variety in general. Moreover, some conjectures
can be verified in low dimensions. There are a lot of papers dedicated to
the classification of various varieties of algebras over the field of
the complex numbers and fields of positive characteristics. However,
 there is also the special field, $\Q_p$, which has zero characteristic.
  There are already works \cite{AyK,FS,Ha,KKO,Ku} devoted to $p$-adic algebras.

In this paper we deal with the low-dimensional Leibniz algebras.
Let us recall that Leibniz algebras present a ``noncommutative'' (to be more precise, a ``nonantisymmetric'')
analogue of Lie algebras \cite{Lo}. These Leibniz algebras satisfy the following
Leibniz identity:
\[ [x,[y,z]]=[[x,y],z] - [[x,z],y] \,.\]

Let $L$ be an arbitrary $n$-dimensional algebra and let $\{e_1, e_2, \dots, e_n\}$
 be a basis of the algebra $L$. Then the table of multiplication on the algebra is defined
by the products of the basic elements, namely,  $[e_i,e_j]=\sum_{k=1}^n \gamma_{i,j}^ke_k$,
where $\gamma_{i,j}^k$ are the structural constants.

From the Leibniz identity, we derive the polynomial equalities for the
structural constants. Thus, the problem of the classification of Leibniz algebras can be
reduced to the problem of the description of the structure constants up to a
nondegenerate basis transformation.

Recall that the description of finite dimensional complex Lie algebras has
been reduced to the classification of nilpotent Lie algebras, which have been completely classified up to dimension 7. In the case of Leibniz algebras the problem of classification of complex Leibniz algebras has been solved till dimension 4.

Note that classifications of 1 and 2-dimensional $p$-adic and complex Leibniz algebras coincide.
 The descriptions of 3-dimensional solvable $p$-adic Leibniz algebras \cite{KKO} and 4-dimensional
  $p$-adic filiform Leibniz algebras \cite{AyK} show that even in these cases the lists of Leibniz
   algebras are essentially wider than in the complex case.
   Moreover, the classification of 5-dimensional $p$-adic filiform Leibniz algebras is known.

 As in the real case to solve a problem in the $p$-adic case there arises
an equation which must be solved in the field of $p$-adic numbers
(see for example \cite{MR1,MR2,VVZ}).  It should be noted that for the
mentioned results it is sufficient to know the solution of the
equations $x^q=a$ for $q=2, 4$. However, in order to classify
6-dimensional filiform $p$-adic Leibniz algebras we also need
the solutions of the equations for $q=3, 5, 6$. In this paper,
using the criteria of solvability for some $p$-adic equations and
structure methods we extend the classification to 6-dimensional $p$-adic Leibniz algebras
and we construct an algorithm running under Maple that it provides criteria of solvability of the  equation
$x^{p^m}=a$, $m \leq p-1$.

\section{Solvability of the equation $x^q=a$}

Let $\Q$ be the field of the rational numbers. Every rational number $x\ne 0$ can be represented
in the form $x = p^r\frac{n}{m}$, where $r, n\in \Z$, $m$ is a positive integer, $(p, n) = 1$, $(p, m) = 1$ and $p$
 is a fixed prime number. The $p$-adic norm of $x$ is given by
\begin{equation*}
|x|_p=
\begin{cases}
p^{-r}, & \text{ for }  x\ne 0;\\
0, & \text{ for }  x = 0.
\end{cases}
\end{equation*}
This norm satisfies the so-called strong triangle inequality
\[|x+y|_p\leq \max\{|x|_p,|y|_p\}\, ,\]
i.e. is a non-Archimedean norm.

The completion of $\Q$ with respect to the $p$-adic norm defines the $p$-adic field
which is denoted by $\Q_p$. Any $p$-adic number $x\ne 0$ can be uniquely represented
in the canonical form
 \[\label{ek}
x = p^{\gamma(x)}(x_0+x_1p+x_2p^2+\cdots),
\]
where $\gamma=\gamma(x)\in \Z$ and $x_j$ are integers, $0\leq x_j \leq p - 1$, $x_0 > 0$, $j = 0, 1,2,\dots $ (see
more detail in \cite{Kh2,VVZ,Sc}). In this case $|x|_p = p^{-\gamma(x)}$.

In this section we consider the equation $x^q=a$ in $\Q_p$, where $p$ is a fixed prime number,
 $q\in \N$ and $0\ne a\in \Q_p$ with canonical form
\[a =p^{\gamma(a)}(a_0 + a_1p +\cdots), \ \ 0\leq a_j \leq p - 1, \ \ a_0 > 0 \,.\]

Solvability of the equation depends on the parameters $p, q, a$.
The following results are known.

\begin{thm}[\cite{Kh2,VVZ}] \label{t1}   The equation
$x^2 = a$ has a solution $x\in \Q_p$, if and only if the following conditions
are fulfilled:
\begin{itemize}
  \item[i)]   $\gamma(a)$ is even;
  \item[ii)] $a_0$ is a quadratic residue modulo $p$ if $p\ne 2$; $a_1 = a_2 = 0$ if $p = 2$.
\end{itemize}
\end{thm}

Let $\eta$  be a unity  (i.e. $|\eta|_p=1$)  which is not a square of any $p$-adic number.

\begin{cor}[\cite{VVZ}] \label{c0} \
\begin{enumerate}
  \item For $p\ne 2$ the numbers
$\varepsilon_1=\eta$, $\varepsilon_2=p$, $\varepsilon_3=p\eta$ are not squares of any $p$-adic number.
\item For $p\ne 2$ any $p$-adic number $x$ can be represented in one of the four following forms:
$x=\varepsilon y^2$, where $y\in \Q_p$ and  $\varepsilon \in \{1,\eta,p,p\eta\}$.
  \item For $p=2$ every $2-$adic number $x$ can be represented in one of the eight following forms:
$x=\varepsilon y^2$, where $y\in \Q_2$ and $\varepsilon\in \{1,2,3,5,6,7,10,14\}$.
\end{enumerate}
\end{cor}

The following theorem gives a generalization of Theorem \ref{t1}.

\begin{thm}[\cite{COR}] \label{t2} Let $q>2$ and $(q,p)=1$. The equation
\[\label{e1}
x^q=a,
\]
has a solution $x\in \Q_p$ if and only if
\begin{enumerate}
  \item $q$ divides $\gamma(a)$;
  \item $a_0$ is a $q$ residue modulo $p$.
\end{enumerate}

\end{thm}

\begin{cor}[\cite{COR}] \label{c1} \
\begin{enumerate}
  \item Let $q$ be a prime number such that $q<p$ and $\eta$ be a unity which is not $q$-th power of some $p$-adic number. Then $p^i\eta^j$, $i,j=0,1,\dots,q-1$ ($i+j\ne 0$), are not $q$-th power of any $p$-adic number.
  \item Let $q$ be a prime number such that $q<p=qk+1$ for some $k\in \N$ and $\eta$ be a unity, which is not $q$-th power of any $p$-adic number. Then any $p$-adic number $x$ can be represented in one of the following form $x=\varepsilon_{ij}y_{ij}^q$ where $\varepsilon_{ij}\in \{p^i\eta^j: i,j=0,1,\dots,q-1\}$ and $y_{ij}\in \Q_p$.
  \item Let $q$ be a prime number such that $q<p\ne qk+1$, for any $k\in \N$. Then any $p$-adic number $x$ can be represented in one of the following forms: $x=\varepsilon_{i}y_{i}^q$, where $\varepsilon_{i}\in \{p^i: i=0,1,\dots,q-1\}$ and $y_{i}\in \Q_p$.
\end{enumerate}
\end{cor}

In the case $q=p$ the solvability condition is given by the following

\begin{thm}[\cite{COR}] \label{t3} Let $q=p$. The equation $x^q=a$,
 has unique solution $x\in \Q_p$ if and only if
\begin{enumerate}
  \item[(i)] $p$ divides $\gamma(a)$;
  \item[(ii)]  $a_0^p \equiv a_0+a_1p \mod p^2$.
\end{enumerate}

 Moreover, the solution can be uniquely obtained from the equations:
 \begin{align} \label{criteria}
x_0= & \ a_0; \notag \\
x_i\equiv & \  a_{i+1}-p^{-1}N_i\mod p, \ \ 0\leq x_i\leq p-1, \ \ i=1,\dots,p-2; \notag \\
x_{pk-1}\equiv & \  a_{pk}-\tilde{N}_{pk}-p^{-1}N_{pk-1}\mod p, \ \ k=1,2,\dots; \notag\\
x_{pk} \equiv  & \ a_{pk+1}+x^{p-2}_0x_1x_{pk-1}\mod p, \ \ k=1,2,\dots;\\
x_{pk+1}\equiv & \ a_{pk+2}-x_0^{p-2}x_1x_{pk-1}-p^{-1}N_{pk+1}\mod p,  \ \ k=1,2,\dots;\notag \\
x_{pk+i}\equiv & \ a_{pk+i+1}-p^{-1}N_{pk+i}\mod p, \ \ i=2,\dots,p-2,\ \ k=1,2,\dots; \notag
\end{align}
where $N_1=0$, $N_2=N_2(x_0,x_1)=\frac{(p-1)p}{2} x_0^{p-2}x_1^2$  and the integer numbers $N_k$  $(k\geq 3)$   depend only on $x_0,x_1,\dots,x_{k-1}$; moreover
$\tilde{N}_{pk}$ depends only on $x_0,x_1,\dots,x_{pk-2}$.
\end{thm}

 \begin{cor}[\cite{COR}] \label{c3} Let $p$ be a prime number.
\begin{itemize}
  \item[(a)] The numbers $\varepsilon \in {\mathcal E}_1=\{1\}\cup\{i+jp: i^p\not\equiv i+jp \mod p^2\}$,  $\delta\in {\mathcal E}_2=\{p^j: j=0,\dots,p-1\}$ and the products $\varepsilon\delta$ are not $p$-th power of any $p$-adic number.
  \item[(b)] Any $p$-adic number $x$ can be represented in one of the following forms
   $x=\varepsilon\delta y^p$, for some $\varepsilon\in {\mathcal E}_1$, $\delta\in {\mathcal E}_2$ and $y\in \Q_p$.
\end{itemize}
  \end{cor}
 In \cite{COR}, we show the following:
 \begin{description}
   \item[$p=3$] Any $x\in \Q_3$ has the form $x=\varepsilon\delta y^3$, where $\varepsilon\in \{1,4,5\}$, $\delta\in \{1,3,9\}$;
   \item[$p=5$] Any $x\in \Q_5$ has the form $x=\varepsilon\delta y^5$, where $\varepsilon\in \{1,11,12,13,14\}$, $\delta\in \{1,5,25,125,625\}$;
   \item[$p=7$] Any $x\in \Q_7$ has the form $x=\varepsilon\delta y^7$, where $\varepsilon\in \{1,8,9,10,11,12,13,\\22,23,24,25,
 26,27,28,36,37,38,39,40,41,42\}$, $\delta\in \{7^i: i=0,1,\dots,6\}$.
 \end{description}

 Now using Theorems \ref{t1}, \ref{t2} and \ref{t3} we shall give solvability conditions of the equation
 $x^q=a$ for $q=mp^s$ for some values $m, s\in \N$ which are necessary for the classification of $p$-adic algebras in the next section.

 \begin{thm}\label{t4}
  Let $q=mp$, $(m,p)=1$. The equation $x^q=a$ has solution $x\in \Q_p$ if and only if
\begin{itemize}
  \item[(i)] $mp$ divides $\gamma(a)$;
  \item[(ii)] $a_0$ is a $m$ residue modulo $p$ and $a_0^p\equiv a_0+a_1p\mod p^2$.
\end{itemize}
\end{thm}
\proof Let us denote $y=x^m$. Then we get the equation $y^p=a$, which by Theorem \ref{t3} has a solution if and only if $p\mid\gamma(a)$ and $a_0^p\equiv a_0+a_1p\mod p^2$. Moreover, $y_0=a_0$ and $\gamma(a)=\gamma(y)p$.  Now, assuming that $y$ is already found, we consider the equation
$x^m=y$, which by Theorem \ref{t2} has a solution if and only if $m\mid\gamma(y)$ and $y_0$ is a $m$ residue $\mod p$. In this case we have $\gamma(y)=m\gamma(x)$.
\endproof

Note that the above mentioned results are sufficient to check of solvability of $x^q=a$ for any $q\in \N$. Indeed, the remaining case is $q=mp^s$ with some $m, s\in \N$, $(m,p)=1$. Denoting $y=x^{p^s}$, we get $y^m=a$ which is the equation considered above. Assume for the last equation the solvability condition is satisfied and its solution is $y={\tilde y}$. Then we have to solve
 $x^{p^s}={\tilde y}$, here we denote $z=x^{p^{s-1}}$ and get $z^p={\tilde y}$. The last equation is the equation considered above. Suppose it has a solution $z={\tilde z}$ then we get $x^{p^{s-1}}={\tilde z}$ which again can be reduced to the previous equation. Iterating the last argument after $(s-1)$ times we obtain the equation $x^p={\tilde a}$ for some ${\tilde a}$.    Consequently,  by this argument we establish the solvability condition of equation $x^{mp^s}=a$ which will be a system of conditions of solvability of equations considered above.

In some particular cases we get the following

\begin{prop} \
\begin{enumerate}
  \item $x^6=a$ has a solution in $\Q_2$ if and only if $6\mid\gamma(a)$ and $a_1=0$;
  \item $x^{(p-1)p}=a$, $p\geq 3$ has a solution in $\Q_p$ if and only if $(p-1)p\mid \gamma(a)$ and $a_0=1$ and $a_1=0$;
  \item $x^{p^2}=a$ has a solution in $\Q_p$ if and only if $p^2\mid \gamma(a)$ and  $a^p_0\equiv a_0+a_1p \mod p^2$ and $a_1=a_2$.
  \item $x^{p^3}=a$ has a solution in $\Q_p$ if and only if  $p^3\mid \gamma(a)$ and  $a^p_0\equiv a_0+a_1p \mod p^2$, $a_1=a_2$ and $a_1\equiv a_3-\frac{p-1}{2}a_0^{p-2}a_1^2\mod p$.
  \item $x^{2^k}=a$ has a solution in $\Q_2$ if and only if  $2^k\mid \gamma(a)$ and  $a_1=\dots=a_{k+1}=0$.
\end{enumerate}
\end{prop}
\proof  One can easily check that for the solvability of an equation $x^q=a$, the
condition $q\mid\gamma(a)$ is necessary for any $q\in \N$.

(1)  This directly follows from Theorem \ref{t4}.

(2)  Using Theorem \ref{t4} we get  $a_0$ is a $p-1$ residue $\mod p$ which is possible only if  $a_0=1$. Then from
 $a_0^p\equiv a_0+a_1p\mod p^2$ we get $a_1=0$.

(3)  Denote $y=x^p$. Then by Theorem \ref{t3}, which provides the solvability of the equation $y^p=a$,  we get $y_0=a_0$, $y_1=a_2$, $a_0^p\equiv a_0+a_1p\mod p^2$ and for $x^p=y$ we have $x_0=y_0$,  $y_0^p\equiv y_0+y_1p\mod p^2$. Consequently $x_0=y_0=a_0$,  $a_0^p\equiv a_0+a_1p\mod p^2$,  $a_0^p\equiv a_0+a_2p\mod p^2$. The last two conditions give $a_1=a_2$.

(4)  Denote  $y=x^{p^2}$. For the equation $y^p=a$ by Theorem \ref{t3} we get $y_0=a_0$, $y_1=a_2$, $y_2\equiv a_3-p^{-1}N_2(y_0,y_1)\mod p$. By part (3) for $x^{p^2}=y$  we have $y_0^p\equiv y_0+y_1p\mod p^2$, $y_1=y_2$. Consequently we get assertion (4).

(5) We use mathematical induction over $k$.  For $k=1$ the assertion is true (see Theorem \ref{t1}). Assume $x^{2^{k-1}}=a$ has solution if and only if $a_1=a_2=\dots=a_k=0$. Denote $y=x^{2^{k-1}}$. Then the equation $x^{2^k}=a$ becomes  $y^2=a$ which has a solution if and only if $a_1=a_2=0$; moreover this solution, $y$, can be obtained by the equations (see \cite{VVZ}):
\begin{equation}\label{p2}\left\{ \begin{array}{ll}
y_0=1 \\
a_3\equiv \frac{y_1(y_1+1)}{2}+y_2 & \mod 2,\\
  y_j=a_{j+1}+N_j(y_0,\dots,y_{j-1}) & \mod 2,  \ \ j\geq 3.
  \end{array} \right.
\end{equation}
By the assumption, the equation $x^{2^{k-1}}=y$ has a solution if and only if $y_1=\dots=y_k=0$.  It is easy to see that $N_j(y_0,\dots, y_{j-1})=N_j(1,0,\dots, 0)=0$, $3\leq j\leq k$, consequently, from the  equations \eqref{p2} we get $a_3=a_4=\dots=a_{k+1}=0$.
\endproof
\begin{cor}\label{c4} \
 \begin{enumerate}
   \item Any $x\in \Q_2$ has the form $x=\varepsilon\delta y^6$, where $\varepsilon\in \{2^i:i=0,1,2,3,4,5\}$, $\delta\in \{1,3\}$, $y\in \Q_2$;
   \item Any $x\in \Q_3$ has the form $x=\varepsilon\delta y^6$, where $\varepsilon\in \{1,2,4,5,7,8\}$, $\delta\in \{3^i: i=0,\dots ,5\}$, $y\in\Q_3$;
   \item Let $\eta$ be a unity which is not a square of any $p$-adic number. Let $\mu$ be a unity which is not a cube of any $p$-adic number.
 Any $x\in \Q_p$, $p\geq 5$, has the form $x=\varepsilon\delta^2y^6$, where $y\in \Q_p$ and $\varepsilon\in \{1,\eta,p,p\eta\}$,
 $\delta\in \{p^i\mu^j: i,j=0,1,2\}$.
 \end{enumerate}
\end{cor}

\begin{rk} \label{rem10} From the results mentioned in this section it follows that for a given $q\in \N$ any $x\in \Q_p$ has
the form $x=\varepsilon y^q$, where $y\in\Q_p$ and $\varepsilon$ varies in a finite set which depends on $p$ and $q$; we denote this set
by ${\mathcal E}_{p,q}$. We have
\[
{\mathcal E}_{p,2}=
\begin{cases}
\{1, \eta, p, p\eta\}, & \text{if $p\ne 2$;}\\
\{1,2,3,5,6,7,10,14\}, & \text{if $p=2$,}
\end{cases}
\]
where $\eta$ is a unity which is not a square of any $p$-adic number.
\[
{\mathcal E}_{p,3}=
\begin{cases}
\{1,2,4\}, & \text{if $p=2$;}\\
\{1,3,4,5,9,12,15,36,45\}, & \text{if $p=3$;}\\
\{p^i\zeta^j: i,j=0,1,2\}, & \text{if $p\geq 5$,}
\end{cases}
\]
where $\zeta$ is a unity which is not a cube of any $p$-adic number.
\[
{\mathcal E}_{p,4}=
\begin{cases}
\big\{\e\delta: \e\in\{1,3,5,7,9,11,13\}, \  \delta\in\{1,2,4,8\}\big\}, & \text{if $p=2$;}\\
\{1,p,p^2,p^3,\eta,p\eta,p^2\eta,p^3\eta\}, & \text{if $p=4k+3, \  k\geq 0$;}\\
\{p^i\eta^j: i,j=0,1,2,3\}, & \text{if $p=4k+1, \ k\geq 1$,}
\end{cases}
\]
where $\eta$ is a unity which is not a square of any $p$-adic number.
\[{\mathcal E}_{p,5}=\left\{\begin{array}{lll}
\{1,p,p^2,p^3,p^4\}, \qquad \qquad  \quad  \quad \mbox{if} \ \ p=2,3, p\ne 5k+1, k\geq 2;\\
\big\{\e\delta: \e\in\{1,11,12,13,14\}, \ \delta\in\{1,5,25,125,625\}\big\},\ \ \mbox{if} \ \ p=5;\\
\{p^i\xi^j: i,j=0,1,2,3,4\}, \quad \quad   \mbox{if} \ \ p=5k+1, \ \ k\geq 2,\\
\end{array}\right.
\]
where $\xi$ is a unity which is not a 5-th power of any $p$-adic number.
\[
{\mathcal E}_{p,6}=
\begin{cases}
\big\{\e\delta: \e\in\{2^i: i=0,1,2,3,4,5\}, \ \delta\in \{1,3\}\big\},  & \text{if $p=2$;}\\
\big\{\varepsilon\delta: \varepsilon\in \{1,2,4,5,7,8\}, \ \delta\in \{3^i: i=0,\dots ,5\}\big\}, & \text{if $p=3$;}\\
\big\{\varepsilon\delta^2: \varepsilon\in \{1,\eta,p,p\eta\}, \ \delta\in \{p^i\mu^j: i,j=0,1,2\}\big\}, & \text{if $p\geq 5$,}
\end{cases}
\]
where $\eta$ is a unity which is not a square of any $p$-adic number and $\mu$ is a unity which is not a cube of some $p$-adic number.
\end{rk}

\section{Computer Program}
In this section we describe a program running under Maple for implementing
the algorithm  \eqref{criteria} discussed in the previous section. The program establishes a criterion of solvability of the  equation
$x^{p^m}=a$, $m \leq p-1$.

\

\noindent \texttt{lis:=proc(p,k)}

\noindent \texttt{\#lis determines the list of summands of the polynomial $N_i$}
\begin{verbatim}
local L, X, Y, Z, i, j, l, ll;
L:=[seq(0,i=0..k-1)]:
X:={}:
Z:={}:
if k= 2 then {[p-2,2]} else
 for i from 0 to iquo(k,k-1) do
   L[k]:=i:
   L[2]:= k - (k-1)*i:
   X:={op(X),L}
end do;
for j from 2 to k-3 do
   for l in X do
       ll:=l;
       for i from 0 to iquo(l[2], k-j) do
           ll[k-j+1]:=i:
           ll[2]:= l[2] - (k-j)*i:
           X:={op(X), ll};
       end do;
    end do;
end do;
for l in X do
        l[1]:=p -`+`(op(2..k,l)):
        Z:={op(Z),l};
end do;
Z;
end if;
end:
\end{verbatim}
\noindent \texttt{pol:=proc(M)}

\noindent \texttt{\# pol determines the polynomial $N_i$ for a list of summands}
\begin{verbatim}
local i,j, X, Y, L;
Y:=0:
for L in M do
  X:=`+`(op(L))!:
  for j from 1 to nops(L) do X:=X/L[j]! end do;
    for i from 0 to nops(L)-1 do
      X:=X*x[i]^(L[i+1])
    od;
  Y:=Y+X;
od;
end:
\end{verbatim}
\verb+pol_N:=(p,k)-> pol(lis(p,k)):+

\noindent\verb+# pol_N+  \texttt{determines the polynomial $N_i$}

\noindent \texttt{sol:=proc(A,p,m)}

\noindent \texttt{\# sol gives the solution of the equation $x^{p^m}=a, \ m \leq p-1$}
\begin{verbatim}
local N, i, j, k;
global L;
N[1]:=(x,y)->0:
k:= nops(A):
for j from 2 to k do
   N[j]:=unapply(pol_N(p,j)/p, seq(x[r],r=0..j)):
end do;
if m = 1 then
   L:=A:             					
       for i from 2 to k-1 do
       L[i]:= A[i+1]-N[i-1](seq(A[j],j=1..i)) mod p
   end do;
   L;
else
   sol(sol(A,p,1),p,m-1)
end if;
end:
\end{verbatim}

\begin{ex}  The criterion of solvability for the equation $x^{p^4}=a$ is:
\begin{itemize}
   \item[(i)] $p^4$ divides $\gamma(a)$;
   \item[(ii)]
\[
\left\{ \begin{array}{lll}
a_0^p \equiv & \  a_0+a_1p & \mod p^2 \\
   a_1 = & \ a_2  & \\
 a_1 \equiv & \  a_3-\frac{p-1}{2} a_0^{p-2}a_1^2  & \mod p^2\\
  a_1 \equiv & \  a_4-\frac{(p-1)(p-2)}{6} a_0^{p-3}a_1^3 +\frac{3(p-1)}{2} a_0^{p-2}a_1^2 & \mod p^2
  \end{array} \right.
\]
\end{itemize}
\end{ex}

\section{Classification of p-adic 6-dimensional filiform Leibniz algebras}

In this section we apply the results of the previous sections. In fact, we use the classification of six-dimensional complex filiform complex Leibniz algebras \cite{Ak} and we omit the cases where complex and $p$-adic considerations lead to the same algebras. Thus, we indicate the cases where equations in complex and $p$-adic numbers are solved in different manner.

Let $L$ be a Leibniz algebra and $L^k$  the lower central series,  $k \geq 1$, defined as follows:
\[L^1=L, \ L^{k+1}=[L^{k},L], \ k\geq 1 \,.\]

\begin{defn} A Leibniz algebra $L$ is said to be filiform if
$dim L^i=n-i$, where $n=dim L$ and $2\leq i \leq n$.
\end{defn}

From the works \cite{GO,OR} we have the description of $n$-dimensional complex filiform Leibniz algebras.
 In fact, the description of such algebras is also true over the field of $p$-adic numbers. Using this fact we conclude that
an arbitrary 6-dimensional $p$-adic filiform Leibniz algebra belongs to one of the following non intersected classes of algebras:
\[
\begin{array}{lll}
\mbox{I.}      &   [e_1,e_1]=e_3, &[e_1,e_2]=\alpha
_1e_4+\alpha_2e_5+\beta e_6, \\
        &[e_i,e_1]=e_{i+1},\ (2\leq i\leq 5),& [e_2,e_2]=
        \alpha_1e_4 + \alpha_2e_5 + \alpha_3e_6, \\
  & & [e_3,e_2] = \alpha_1e_5 + \alpha_2e_6, \\
  &&  [e_4,e_2] = \alpha_1e_6,\\[3mm]
\mbox{II.} &[e_1,e_1] = e_3,& [e_1,e_2] = \beta_1e_4 + \beta_2e_5 + \beta_3e_6, \\
  & [e_i,e_1]= e_{i + 1},\ (3\leq i\leq 5),& [e_2,e_2]=\gamma e_6, \\
    & &   [e_3,e_2] = \beta_1e_5 + \beta_2e_6,\\
    & &   [e_4,e_2] = \beta_1e_6,\\[3mm]
 \mbox{III.}& [e_i,e_1]=e_{i+1}, \ (2\leq i \leq 5), & [e_1,e_i]=-e_{i+1}, \ (3\leq i \leq 5), \\
&[e_1, e_1]=\theta_1e_6, & [e_1, e_2]=-e_3+\theta_2e_6,\\
&[e_2,e_2]=\theta_3e_6, & [e_2,e_3]=-[e_3,e_2]=\alpha e_5+\beta e_6,\\
&[e_2,e_4]=-[e_4,e_2]=\alpha e_6, &[e_3, e_4]=-[e_4, e_3]=\delta e_6,\\
&[e_2, e_5]=-[e_5,e_2]=-\delta e_6, & \\
\end{array}
\]
where $\delta \in \{0, 1\}$, $\{e_1, e_2, e_3, e_4, e_5, e_6\}$ is a basis of the algebra and omitted products are equal to zero.

We denote by $L_1(\alpha_1,\alpha_2,\alpha_3,\beta)$, $L_2(\beta_1,\beta_2,\beta_3,\gamma)$ and $L_3(\t_1,\t_2,\t_3,\alpha,\beta,\delta)$ the algebras from the classes I, II and III, respectively.

 Below, we consider the class III (which contains Lie algebras) since several equations,
  which were considered in the previous section (i.e., the solutions of which are different in complex and $p$-adic numbers) appear there.
  Moreover, in this class, the equations of the form $x^q=a$ only appear, where $q$ is a product of prime numbers.

 Let us take the general transformation of the generator basis elements of the algebra $L_3(\t_1,\t_2,\t_3,\alpha,\beta,\delta)$  of the class III in the following form:
\begin{align*}
e_1'= & \ A_1e_1+A_2e_2+A_3e_3+A_4e_4+A_5e_5+A_6e_6,\\
 e_2'= & \ B_1e_1+B_2e_2+B_3e_3+B_4e_4+B_5e_5+B_6e_6,
\end{align*}
where $A_1B_2-A_2B_1\neq 0$.

Then we generate the rest of the elements of the basis $\{e_3', e_4', e_5', e_6'\}$ and the parameters of the algebra $L_3(\theta_1', \theta_2', \theta_3', \alpha', \beta', \delta')$ in the new basis are expressed via the parameters $\{\theta_1, \theta_2, \theta_3, \alpha, \beta, \delta\}$ in the following equalities:
\begin{gather*}
\theta_1'=\frac{A_1^2\theta_1+A_1A_2\theta_2+A_2^2\theta_3}{A_1^3(A_1+A_2\delta)B_2}, \
\theta_2'=\frac{A_1\theta_2+2A_2\theta_3}{A_1^3(A_1+A_2\delta)}, \
\theta_3'=\frac{B_2\theta_3}{A_1^3(A_1+A_2\delta)}\, ,\\
\alpha'=\frac{B_2}{A_1^2}\alpha, \
\beta'=\frac{\beta A_1^2B_2^2+2\alpha^2A_1A_2B_2^2+\alpha^2\delta A_2^2B_2^2+\delta A_1^2B_3^2-2\delta A_1^2B_2B_4}{A_1^4B_2(A_1+A_2 \delta)} \, ,\\
\delta'=\frac{B_2 \delta}{A_1+A_2\delta} \, ,
\end{gather*}
where $A_1B_2\neq0$.
The next step in the classification is the analysis of the possible cases.

Below, we consider only cases in which there appears a difference between complex and $p$-adic solutions of the equations.
\begin{description}
  \item[1. Case $\delta=\alpha=\beta=\theta_3=0,  \ \theta_2\ne 0$] The above equalities have the form $\theta_1'=\frac{A_1\theta_1+A_2\theta_2}{A_1^3B_2}, \ \theta_2'=\frac{\theta_2}{A_1^3}$.
Let $\theta_2=\varepsilon y^3$, where $\varepsilon \in {\mathcal E}_{p,3}$. Putting $A_1=y$ and $A_2=-\frac{\theta_1}{\varepsilon y^2}$ we obtain $\theta_2'=\varepsilon, \ \theta_1'=0$, and so we get the algebra $L_3(0,\varepsilon,0,0,0,0)$, $\varepsilon \in {\mathcal E}_{p,3}$.
  \item[2. Case $\delta=\alpha=\beta=0, \ \theta_3(\theta_2^2-4\theta_1\theta_3)\neq0$] Then similarly as in the complex case we derive
$\theta_2'=0, \ \theta_3'=1$ and $\theta_1'=\frac{4\theta_1\theta_3-\theta_2^2}{4A_1^6}$.
Let $\frac{\theta_1\theta_3-\frac{\theta_2^2}{4}}{A_1^6}= \varepsilon y^6$. Taking $A_1=y$, we obtain $\theta_1'=\varepsilon$ and the algebra $L_3(\varepsilon,0,1,0,0,0)$, $\varepsilon \in {\mathcal E}_{p,6}$.
\item[3. $\delta=\alpha=\theta_2=\theta_3=0, \ \beta\theta_1\neq 0$] Then $\beta'=1,\  \theta_3'=0$ and $p$-adic property appears in the equation $\theta_1'=\frac{\theta_1\beta}{A_1^5}$. If $\theta_1\b=\varepsilon y^5$, then putting
 $A_1=y$ we get $\theta_1'=\varepsilon$ and the  algebra $L_3(\varepsilon,0,0,0,1,0)$, $\varepsilon \in {\mathcal E}_{p,5}$.

\item[4. Case $\delta=\alpha=\theta_3=0, \  \beta\theta_2\neq 0$] Then the above expressions reduce to the following: $\theta_1'=\frac{(A_1\theta_1+A_2\theta_2)\beta}{A_1^6}, \ \theta_2'=\frac{\theta_2}{A_1^3}$. If $\theta_2=\varepsilon y^3$, then by making $A_1=y, \ A_2=-\frac{\theta_1}{\varepsilon y^2}$ we obtain $\theta_1'=0, \ \theta_2'=\varepsilon$ and the algebra $L_3(0,\varepsilon,0,0,1,0)$, $\varepsilon \in {\mathcal E}_{p,3}$.
\item[5. Case $\delta=\theta_2=\theta_3=0, \ \alpha\ne 0$] Then we deduce
$\alpha'=1, \ \theta'_2=\theta'_3=0$ and $\theta'_1=\frac{\theta_1}{A_1^4}\alpha$, $\beta'=\frac{A_1\beta+2\alpha^2A_2}{A_1^2\alpha}$. Let $\theta_1\alpha=\varepsilon y^4$. Put $A_1=y, \ A_2=-\frac{y\beta}{2\alpha^2}\  \Rightarrow \ \theta_1'=\varepsilon \beta'=0$. Thus, we get the algebra $L_3(\varepsilon,0,0,1,0,0)$, $\varepsilon\in {\mathcal E}_{p,4}$.
\item[6. Case $\delta=\theta_3=0, \ \alpha\theta_2\neq0$] In this case we deduce $\alpha'=1$, $\theta'_1=\frac{A_1\theta_1+A_2\theta_2}{A_1^5}\alpha$, $\theta'_2=\frac{\theta_2}{A_1^3}$ and $\beta'=\frac{A_1\beta+2A_2\alpha^2}{A_1^2\alpha}$. If $\theta_2=\varepsilon y^3$, then making $A_1=y, \ A_2=-\frac{y\beta}{2\alpha^2}$ we derive $\theta_2'=\varepsilon, \ \beta'=0$ and $\theta_1'=\frac{2\alpha^2\theta_1-\beta\theta_2}{2\alpha^2y^4}$ as parameter, i.e. the algebra $L_3(\theta_1',\varepsilon,0,1,0,0)$, $\varepsilon\in {\mathcal E}_{p,3}$ is obtained.
\item[7. Case $\delta=0,  \ \alpha\theta_3\neq0$] Then we derive $\alpha'=1$ and the mentioned equalities have the form:
$\theta_1'=\frac{(A_1^2\theta_1+A_1A_2\theta_2+A_2^2\theta_3)\alpha}{A_1^6}$, $\theta_2'=\frac{A_1\theta_2+2A_2\theta_3}{A_1^4}$, $\theta'_3=\frac{\theta_3}{A_1^2\alpha}$ and  $\beta'=\frac{A_1\beta+2A_2\alpha^2}{A_1^2\alpha}$. Assume $\frac{\theta_3}{\alpha}=\varepsilon y^2$ and taking $A_1=y, \ A_2=-\frac{y\beta}{2\alpha^2}$ then  $\beta'=0, \ \theta_3'=\varepsilon$ and $\theta_1', \ \theta_2'$ are parameters. Consequently, we have the  algebra $L_3(\theta_1',\theta_2',\varepsilon,1,0,0), \ \varepsilon\in{\mathcal E}_{p,2}$.
\item[8. Case $\delta=1, \theta_2=\theta_3=0, \ \alpha\theta_1\neq0$]
 Then we get $\theta_2'=\theta_3'=\beta'=0$ and the mentioned equalities have the form:
$\theta_1'=\frac{\theta_1}{A_1(A_1+A_2)^2}, \ \alpha'=\frac{A_1+A_2}{A_1^2}\alpha$.
Let $\alpha^2\theta=\varepsilon y^5$, then assuming $A_1=y,
\ A_2=\frac{A_1(A_1-\alpha)}{\alpha}$, we obtain the algebra $L_3(\varepsilon, 0, 0, 1, 0, 1)$ for $\varepsilon\in{\mathcal E}_{p,5}$.
\item[9. Case $\delta=1, \theta_3=0, \ \theta_2(\theta_1-\theta_2)\neq0$]
Similarly to complex case we get $\delta'=\theta_2'=1, \ \beta'=\theta_3'=0$ and $\theta_1'=\frac{A_1^3(\theta_1-\theta_2)+\theta_2^2}{\theta_2^2}, \ \alpha'=\frac{\theta_2}{A_1^4}\alpha$.

If $\alpha\neq0$, then assuming $\alpha\theta_2=\varepsilon y^4$ and $A_1=y$,
 we obtain $\alpha'=\varepsilon, \ \theta_1'$ parameter and the algebra $L_3(\theta_1', 1, 0,
  \varepsilon, 0, 1), \ \varepsilon\in{\mathcal E}_{p,4}$.

If $\alpha=0$, then assuming $\frac{\theta_2^2}{\theta_2-\theta_1}=\varepsilon y^3$ and $A_1=y$,
we obtain $\alpha'=0, \ \theta_1'=1-\frac{1}{\varepsilon}$ and the algebra $L_3(1-\frac{1}{\varepsilon}, 1, 0, 0, 0, 1), \ \varepsilon\in{\mathcal E}_{p,3}$.
\item[10. Case $\delta=1, \ \theta_3(2\theta_3-\theta_2)\neq0$]
Suppose that $\theta_3=\varepsilon y^3$. Putting $A_1=y, \ A_2=-\frac{\theta_2}{2\varepsilon y^2}$, we get $\theta_2'=0, \ \theta_3'=\varepsilon$ and $\alpha', \theta_1'$ are parameters. So, in this case, we get the algebra $L_3(\theta_1', 0, \varepsilon, \alpha', 0, 1), \ \varepsilon\in{\mathcal E}_{p,3}$.
\item[11. Case $\delta=1,  \ 2\theta_3=\theta_2, \theta_3\neq0$]
Suppose that $\theta_3=\varepsilon y^3$, then $\theta_2=2\varepsilon y^3$. Putting $A_1=y$, we get $\theta_2'=2\varepsilon$ and $\theta_3'=\varepsilon$.

If $4\theta_1\theta_3-\theta_2^2=0$ then $\theta_1'=\e$ and $\alpha'=\frac{(y+A_2)\alpha}{y^2}$.
In the case when $\alpha=0$ we get the algebra $L_3(\e,2\e,\e,0,0,1)$, $\e\in {\mathcal E}_{p,3}$.

In the opposite case, i.e., $\alpha\ne 0$ we put $A_2=\frac{y(y-\alpha)}{\alpha}$ which implies  $\alpha'=1$. Thus the algebra $L_3(\e,2\e,\e,1,0,1)$ is obtained.

If $4\theta_1\theta_3-\theta_2^2\ne 0$ then assuming $4\theta_1\theta_3-\theta_2^2=\nu z^2$ and $A_2=\frac{z}{y^2}-y$ we deduce  $\theta_1'=\frac{\nu-4\e^2}{4\e}$ and $\alpha'$ is a parameter.
In this subcase we obtain  the algebra $L_3\big(\frac{\nu-4\e^2}{4\e},2\e,\e,\alpha',0,1\big)$, $\e\in {\mathcal E}_{p,3}$, $\nu\in {\mathcal E}_{p,2}$.
\end{description}
\begin{rk} Note that if $\e=\e'y^q$ for some $y\in \Q_p$, then the algebras, which were obtained
 in the case where  similar relations between new and old parameters appear, are isomorphic.
  Take into account that $12=11y^5$,  $13=11z^5$,  $14=11t^5$,  for some $y,z,t\in \Q_5$
  and $4=5y^3$,  $12=5z^3$, $36=5t^3$ for some $y,z,t\in \Q_3$ in the classification list
   we will omit $4, 12,36$ from $ {\mathcal E}_{3,3}$ and $12, 13, 14$ from ${\mathcal E}_{5,5}$.
\end{rk}
Let us denote
 \[\tilde{{\mathcal E}}_{3,3}=\{1,3,5,9,15,45\}, \quad \tilde{{\mathcal E}}_{5,5}=\big\{\e\delta: \e\in\{1,11\},  \ \ \delta\in\{1,5,25,125,625\}\big\}\, . \]
Applying similar $p$-adic argumentations to the classes I, II and adding algebras,
 which are obtained as for the case of the complex field we present the following classification list in the next statement.

\begin{thm}\label{t20} An arbitrary $p$-adic six-dimensional filiform Leibniz algebra is isomorphic to one of the
following pairwise non-isomorphic algebras:

\

\begin{tabular}{|p{3cm}|p{3.08cm}|p{3.8cm}|p{2.9cm}|}
                 \hline
               $L_1(0,0,0,0)$ & $L_1(0,0,1,1)$ & $L_1(0,1,0,0)$ & $L_1(0,1,1,\beta)$  \\ \hline
                 $L_1(1,0,\alpha,\beta)$ & $L_1(1,-2,5,5)$ & $L_1(0,0,0,\e)$ \quad \quad  \  \ $\e\in {\mathcal E}_{p,3}$, $p\ne 3$, $\e\in\tilde{{\mathcal E}}_{3,3}$ & $ L_1(0,0,1,\e+1)$  $\e\in {\mathcal E}_{p,3}, p\ne 3$, $\e\in\tilde{{\mathcal E}}_{3,3}$\\\hline
                 $L_1(1,-2,5,\e+5)$ \ $\e\in {\mathcal E}_{p,3}, p\ne 3$, $\e\in\tilde{{\mathcal E}}_{3,3}$ & $ L_1(1,-2,\nu+5,\beta)$ $\nu\in {\mathcal E}_{p,2}$ &  & \\\hline
                 $ L_2(0,0,0,0)$ & $L_2(0,0,1,0)$ & $L_2(0,1,0,0)$ & $L_2(0,1,1,0)$ \\\hline
                  $ L_2(0,1,0,1)$ & $L_2(1,0,0,0)$ & $L_2(1,0,\beta,1)$ & $L_2(1,0,\nu,0)$  \  \ $\nu\in {\mathcal E}_{p,2}$ \\\hline
                   $L_3(0,0,0,0,0,0)$ & $ L_3(1,0,0,0,0,0)$ & $L_3(0,0,1,0,0,0)$ & $L_3(0,0,0,0,1,0)$  \\\hline $L_3(\alpha,0,1,0,1,0)$ &
                      $ L_3(0,0,0,1,0,0)$ & $L_3(0,0,0,0,0,1)$ & $L_3(0,0,0,1,0,1)$  \\\hline
                      $L_3(1,0,0,0,0,1)$ &
                        $ L_3(1,1,0,0,0,1)$ & $L_3(\alpha,\beta,\e,1,0,0)$  \ $\e \in {\mathcal E}_{p,3}$
                         & $L_3(0,\e,0,0,0,0)$ $\e \in {\mathcal E}_{p,3}$ \\\hline
                          $L_3(0,\e,0,0,1,0)$  $\e \in {\mathcal E}_{p,3}$ &
                          $ L_3(\alpha,\e,0,1,0,0)$  $\e \in {\mathcal E}_{p,3}$ & $L_3\big(\frac{\e-1}{\e},1,0,0,0,1\big)$
                           $\e \in {\mathcal E}_{p,3}$
                          & $ L_3(\alpha,0,\e,\beta,0,1)$  $\e \in {\mathcal E}_{p,3}$  \\\hline
                          $L_3(\e,2\e,\e,0,0,1)$ $\e \in {\mathcal E}_{p,3}$ & $ L_3(\e,2\e,\e,1,0,1)$ $\e \in {\mathcal E}_{p,3}$ & $L_3\big(\frac{\nu-4\e^2}{4\e},2\e,\e,\alpha,0,1\big)$ $ \e\in {\mathcal E}_{p,3}$, $\nu\in {\mathcal E}_{p,2}$  & $ L_3(\xi,0,0,1,0,0)$  $\xi\in {\mathcal E}_{p,4}$ \\\hline
               \end{tabular}

\begin{tabular}{|p{3cm}|p{3.08cm}|p{3.8cm}|p{2.9cm}|}
                 \hline
                              $ L_3(\alpha, 1, 0,\xi,0,1)$ $\xi\in {\mathcal E}_{p,4}$ & $L_3(\zeta,0,0,0,1,0)$ $\zeta\in {\mathcal E}_{p,5}$,
  $p\ne 5$, $\zeta\in \tilde{{\mathcal E}}_{5,5}$ & $  L_3(\zeta,0,0,1,0,1)$ \ \ $\zeta\in {\mathcal E}_{p,5}$,
  $p\ne 5$, $\zeta\in \tilde{{\mathcal E}}_{5,5}$ & $ L_3(\mu,0,1,0,0,0)$ $\mu\in {\mathcal E}_{p,6}$ \\\hline
               \end{tabular}
\

where $\alpha, \beta\in \Q_p$.
\end{thm}
\section*{Acknowledgements}

We would like to thank Felipe Gago for his useful help during the implementation of the program.

{}

 \address{\small \rm  Manuel Ladra: Departament of Algebra,  University of Santiago de Compostela, 15782
Santiago de Compostela, Spain}\\ \email{manuel.ladra@usc.es}

 \address{\small \rm Bakhrom A. Omirov: Institute of Mathematics and Information Technologies,
Tashkent, Uzbekistan} \\ \email {omirovb@mail.ru}

 \address{\small \rm Utkir A. Rozikov: Institute of Mathematics and Information Technologies,
Tashkent, Uzbekistan} \\ \email {rozikovu@yandex.ru}
\end{document}